\RequirePackage{fix-cm}
\documentclass[twocolumn]{svjour3}          
\smartqed  

\usepackage{graphicx}
\usepackage{amsmath,amsfonts,amssymb}
\usepackage{graphicx}
\usepackage{wrapfig}
\usepackage{url}
\usepackage{multirow}
\usepackage[titletoc]{appendix}
\usepackage{lipsum}
\usepackage{mdframed}
\usepackage{fancyvrb}
\usepackage{tabu}
\usepackage{tabularx}
\usepackage{booktabs}
\usepackage{multirow}
\usepackage{xspace}
\usepackage{enumitem}
\usepackage{algorithm}
\usepackage{algorithmic}
\usepackage[caption=false]{subfig}
\usepackage{hyperref}

\usepackage{dsfont}

\newcommand{\R}{\mathds{R}}
\usepackage[T1]{fontenc}
\usepackage[utf8]{inputenc}
\usepackage{siunitx}

\usepackage{todonotes}

\newcommand{\divergence}{\mathop\mathrm{div}}

\newcommand{\Pol}{\mathbb{P}}
\renewcommand{\R}{\mathbb{R}}
\newcommand{\dive}{\mathop\mathrm{div}}
\newcommand{\TOL}{\mathrm{TOL}}

\journalname{International Journal of Computing and Visualization in Science and Engineering, 31, 12, 2022}

\begin{document}

\title{Efficient adaptivity for simulating cardiac electrophysiology with spectral deferred correction methods\thanks{This work was supported by the European High-Perfor\-mance Computing Joint Undertaking EuroHPC under grant agreement No 955495 (MICROCARD) co-funded by the Horizon 2020 programme of the European Union (EU) and the German Federal Ministry of Education and Research. }
}


\author{Fatemeh Chegini \and Thomas Steinke \and Martin Weiser}



\institute{Fatemeh Chegini \at
           Zuse Institute Berlin\\
           \email{chegini@zib.de}
           \and
           Thomas Steinke \at
           Zuse Institute Berlin\\
           \email{steinke@zib.de}
           \and
           Martin Weiser \at
           Zuse Institute Berlin\\
           \email{weiser@zib.de}
}

\date{}

\maketitle

\begin{abstract}
The locality of solution features in cardiac electrophysiology simulations calls for adaptive methods. Due to the overhead incurred by established mesh refinement and coarsening, however, such approaches failed in accelerating the computations. Here we investigate a different route to spatial adaptivity that is based on nested subset selection for algebraic degrees of freedom in spectral deferred correction methods. This combination of algebraic adaptivity and iterative solvers for higher order collocation time stepping realizes a multirate integration with minimal overhead. This leads to moderate but significant speedups in both monodomain and cell-by-cell models of cardiac excitation, as demonstrated at four numerical examples.


\end{abstract}

\section{Introduction}
\label{sec:introduction}

Cardiac arrhythmia causes about $15\%$ of all mortality, primarily due to disorganization of cardiac tissue at the cellular scale, impacting the propagation of myocardium excitation. In order to understand disease mechanisms, diagnose illness, and design effective drugs for treatments, numerical modeling of cardiac electrophysiology is necessary. Since different physiological mechanisms of interest, several models of different complexity and level of detail are in use~\cite{ColPavSca2014}. 

The coarsest description of cardiac excitation is provided by eikonal models describing the activation time directly~\cite{doi:10.1137/S0036139901389513}. They are efficient to solve but provide only the activation patterns. More detailed information is available in the homogenized monodomain model describing the evolution of the transmembrane voltage, ion concentrations, and states of ion channels~\cite{leon1991computer}. In particular for strongly differing extracellular and intracellular conductivities, the bidomain model~\cite{tung1978bi,colli1990wavefront}, treating the extracellular potential explicitly, achieves even higher fidelity at a further increased computational effort, and is solved routinely on compute clusters and GPUs~\cite{neic2012accelerating}.

These models describe the myocardium as a single homogeneous medium, or an overlay of two such media in case of the bidomain model, and therefore cannot capture effects of myocyte size and shape, cellular inhomogeneity and connectivity patterns, or spatially varying ion channel density. Such tissue properties require a cellular resolution to be faithfully represented. Consequently, heterogeneous models on the cellular scale such as the EMI model (extracellular-membrane-intracellular) have been proposed and investigated recently~\cite{becue:hal-01910679,tveito2021tris,JaegerTveito2021}. With the increased spatial and temporal resolution comes a tremendous increase in computational effort for simulation of the cardiac excitation, calling on one hand for high performance computing facilities, and on the other hand for more efficient algorithms.

The solutions of mono- and bidomain as well as EMI models exhibit traveling depolarization and repolarization fronts, which are thin compared to organ scale. This locality of solution features makes adaptive spatial discretizations attractive in view of reducing computational effort. Classical spatio-temporal mesh adaptivity in time stepping schemes has been proposed~\cite{doi:10.1137/050634785,bendahmane2010multiresolution,BELHAMADIA2022101656} and reported to reduce the discretization size in terms of number of degrees of freedom (dofs) by a large factor, but found to be ineffective in reducing the overall computational effort. This is due to the overhead incurred by error estimation, frequent mesh refinement and coarsening, as well as repeated assembly of mass and stiffness matrices. Block-based multiresolution sche\-mes~\cite{KRAUSE201579} achieve a higher efficiency at the expense of the discretization being less precisely adapted to the actual solution. Consequently, method of lines approaches with fixed spatial discretization are ubiquitous despite their wastefully fine discretization in the bulk of the domain.

In this paper, we investigate a different approach to spatio-temporal adaptivity, first proposed in~\cite{WeiserChegini2022}, that is executed completely on the algebraic level and makes use of extremely fast and simple a posteriori error estimates for dof selection. It relies on the combination of spatial adaptivity with higher order time integration with spectral deferred correction (SDC) methods~\cite{dutt2000spectral,minion2003semi,Weiser2015}. Those are stationary iterative solvers for collocation systems, and offer a high flexibility for combination with adaptivity and inexact computation~\cite{SpRuMiEmKr2016,WeiserGhosh2018,weiser2014spectral}. In cardiac electrophysiology we observe that the significant support of SDC corrections shrinks over the iteration. By interleaving the SDC iteration with a progressive spatial subdomain restriction we reduce the computational effort of later iterations while respecting an overall requested tolerance. This cheap algebraic adaptivity by local SDC truncation can also be interpreted as a kind of natural predictor-corrector multirate integration in the spirit of~\cite{SavcencoHundsdorferVerwer2007}, but gains efficiency from interleaving with the SDC iteration.

The remainder of the paper is organized as follows. In Sec.~\ref{sec:main} we define the EMI model and the monodomain models and describe their spatial discretization with finite elements resulting in a large scale ordinary differential equation (ODE). Sec.~\ref{sec:SDC} is devoted to the time discretization of the ODE with implicit-explicit operator splitting as basic method and the SDC iteration resulting in a higher order method. In Sec.~\ref{sec:adaptivity}, the algebraic adaptivity concept is introduced and theoretically justified. A cheap a posteriri error estimator based on the linear convergence of SDC methods is worked out and completes the adaptive scheme. Numerical examples for monodomain and EMI problems are given in Sec.~\ref{sec:experiments}, demonstrating the improved performance, and discussed in Sec.~\ref{sec:discussion}.





\section{Problem definition}
\label{sec:main}

Excitation propagation in the myocardium is facilitated by diffusion of ions in the intracellular and the extracellular space, and their transport across the cell membranes through various ion channels with specific nonlinear dynamics. We will first describe a detailed model capturing the cellular geometry and then move to coarser models that can bederived by mathematical homogenization.

\subsection{The EMI model}

\begin{figure}
  \centering
  \includegraphics[width=0.8\linewidth]{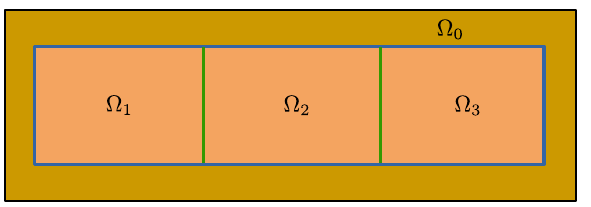}
  \caption{A simple EMI domain comprising three myocytes  $\Omega_i$, $i=1,\dots,3$, and the extracellular domain $\Omega_0$. The green lines represent gap junction interfaces connecting different myocytes, the blue lines represent membranes between myocytes and extracellular medium.}
  \label{fig:cells}
\end{figure}

The EMI (extracellular-membrane-intracellular) model of electrophysiology~\cite{becue:hal-01910679,tveito2021tris,JaegerTveito2021} describes the myocardium as a collection of pairwise disjoint myocytes $(\Omega_i)_{i=1,\dots,N}$ which, together with the extracellular space $\Omega_0$, cover the whole domain $\Omega\subset\R^d$, $d\in\{2,3\}$, occupied by the myocardium, i.e. $\overline\Omega = \bigcup_{i=0}^N \overline\Omega_i$, see Fig.~\ref{fig:cells} for a sketch. 
Ions can diffuse within each myocyte and in the extracellular domain, subject to conductivities $\sigma_i$, which leads to electric intra- and extracellular potentials $u_i\in H^1(\Omega_i)$. Ion currents across the membranes are due to passive ion channels in the case of gap junctions between adjacent myocytes or due to active or passive ion channels controlling the exchange of ions between myocytes and the extracellular space. The transmembrane current $n^T\sigma_i\nabla u_i$  consists of the ion current $I_{ij}^{\rm ion}$ and the capacitive current $C_m\dot v_{ij}$. The ion current depends on the transmembrane voltage $v_{ij} = u_i - u_j$ defined on the membrane $F_{ij} = \partial\Omega_i \cap \partial\Omega_j$ separating adjacent subdomains $\Omega_i$ and $\Omega_j$ as well as the state $w_{ij}$ of ion channels, which in turn follows a nonlinear dynamic given by $R$. We consider linear gap junctions here, i.e. the ion channel state $w_{ij}$ is only effective between intra- and extracellular space ($ij=0$).
This setting results in the partial differential algebraic system
\begin{equation} \label{eq: cell-by-cell model}
		\begin{cases}
			- \dive (\sigma_i \nabla u_i) = 0		&\quad \text{in }\Omega_i\\
			- n_i^T \sigma_i \nabla u_i = C_m \dot v_{ij}  + I_{ij}^{\rm ion}(v_{ij},w_{ij})			&\quad \text{on $ F_{ij}, \; i\ne j	$}	\\
            n^T\sigma_i\nabla u_i  +\epsilon u_i = 0 & \quad \text{on $\partial\Omega_i \cap \partial\Omega$} \\
			\dot w_{ij} = R(v_{ij}, w_{ij}) & \quad \text{on $F_{ij}, ij=0$} \\
            w_{ij} = 0 & \qquad ij > 0.
		\end{cases}
\end{equation}	
Here, $n_i$ denotes the unit outer normal of $\Omega_i$ where $i=0,\dots, N$. Note that $v_{ij} = -v_{ji}$, $I_{ij}^{\rm ion}(v_{ij},w_{ij}) = -I_{ji}^{\rm ion}(v_{ji},w_{ji})$, but $w_{ij} = w_{ji}$. The small value $\epsilon>0$ in the Robin boundary condition on $\partial\Omega$ makes the solution unique, which for a pure Neumann problem would be defined only up to a constant.

\paragraph{Weak formulation}

Multiplying~\eqref{eq: cell-by-cell model} by a test function $\phi_i\in H^1(\Omega_i)$ and integrating by parts yields
\begin{align*}
		0 &= - \int_{\Omega_i} \dive (\sigma_i \nabla u_i) \, \phi_i \, dx	\\	
		&= \int_{\Omega_i} \nabla u_i^T \sigma_i  \nabla \phi_i  \, dx 	- 	\int_{\partial \Omega_i} n_i^T \sigma_i \nabla u_i \phi_i \, ds	\\
		&= \int_{\Omega_i} \nabla u_i^T \sigma_i \nabla \phi_i \, dx 	\\
  &+ \sum_{i\ne j} \int_{F_{ij}} \left( C_m \dot v_{ij} + I_{ij}^{\rm ion}(v_{ij},w_{ij}) \right) \phi_i \, ds
			\quad
		\forall \phi_i \in H^1(\Omega_i).
\end{align*}
Summing over all subdomains yields the weak formulation for $u\in V=\prod_{i=0}^N H^1(\Omega_i)$,
\begin{align}\label{eq:weak-EMI}
    0 &=\sum_{i=0}^N \int_{\Omega_i} \nabla u_i^T \sigma_i \nabla \phi_i \, dx \\ 
    &+ \sum_{i=0}^N \sum_{i\ne j} \int_{F_{ij}} \left( C_m (\dot u_i-\dot u_j) + I_{ij}^{\rm ion}(u_i-u_j,w_{ij}) \right) \phi_i \, ds 
\end{align}
for all $\phi\in V$.

For brevity, we define the symmetric bilinear forms
\begin{align*}
    a(u,\phi) = \sum_{i=0}^N \int_{\Omega_i} \nabla u_i^T \sigma_i \nabla \phi_i \, dx , \\
    \quad
    p(\dot u,\phi) = \sum_{i=0}^N \sum_{i\ne j} \int_{F_{ij}} C_m (\dot u_i-\dot u_j) \phi_i \, ds , 
\end{align*}
and the linear form
\[
    f(\phi;u,w) = \sum_{i=0}^N \sum_{i\ne j} \int_{F_{ij}}  I_{ij}^{\rm ion}(u_i-u_j,w_{ij}) \phi_i \, ds,
\]
such that we can write~\eqref{eq:weak-EMI} as
\[
    a(u,\phi)+p(\dot u,\phi) + f(\phi;u,w) = 0 \quad \forall\phi\in V.
\]

\paragraph{Finite element discretization}
Let $\mathcal{T}$ be a conforming simplicial triangulation of $\Omega$ compatible with the subdomain structure, i.e. each $T\in\mathcal{T}$ is contained in the closure of exactly one subdomain $\Omega_i$. We consider the standard finite element spaces $V_{h,i}=\{u\in C(\overline\Omega_i) \mid \forall T\in\mathcal{T}, T\subset\overline\Omega: u|_T\in\Pol_p\}$ of piecewise polynomial and globally continuous functions and the product space $V_h = \prod_{i=0}^N V_{h,i}$ containing finite element functions which are discontinuous across the cell membranes. Using a Lagrange basis $(\varphi_k)_{k=1,\dots,N_h}$ of $V_h$, the standard Galerkin approach 
\[
    a(u_h,\varphi_k) + p(\dot u_h,\varphi_k) + f(\varphi_k;u_h,w_h) = 0, \quad k=1,\dots,N_h
\]
transforms the weak formulation~\eqref{eq:weak-EMI} into a large scale ordinary differential-algebraic equation
\begin{equation} \label{eq: abstract reaction diffusion}
    M\dot u_h = -A u_h - b(u_h,w_h)
\end{equation}
of index one, with $M_{kl} = p(\varphi_k,\varphi_l)$, $A_{kl}=a(\varphi_k,\varphi_l)$, $b(u_h,w_h) = f(\varphi_k;u_h,w_h)$, and the slight abuse of notation of identifying $u_h\in V_h$ with its coefficient vector $u_h\in\R^{N_h}$ with respect to the basis $(\varphi_k)_k$. On the cell membranes between myocytes and extracellular domain only, the pointwise ordinary differential equations (ODE) for $w_{ij}$ yield the corresponding ODEs 
\begin{equation} \label{eq: abstract gating}
    \dot w_h = R(u_h,w_h)
\end{equation}
for the finite element coefficients $w_h\in\R^{n_h}$. Note that the dimension $N_h$ of~\eqref{eq: abstract reaction diffusion} is usually much larger than the dimension $n_h$ of~\eqref{eq: abstract gating}, since the gating variables $w_h$ are restricted to the nodes on the cell membranes.

\subsection{The bidomain and monodomain models}

Assuming a periodic arrangement of myocytes, a mathematical homogenization approach~\cite{tung1978bi} yields the bidomain system for the extracellular potential $u_{\rm ext}$ and the intracellular potential $u_{\rm int}$. In the homogenized domain $\Omega$, where at every spatial point both intracellular and extracellular space coexist, the potentials satisfy the system
\begin{align*}
    \beta\left(C_m \dot v + I^{\rm ion}(v,w)\right) &= \divergence(\sigma_{\rm int} \nabla u_{\rm int}) \\
    -\beta\left(C_m \dot v + I^{\rm ion}(v,w)\right) &= \divergence(\sigma_{\rm ext} \nabla u_{\rm ext}) \\
    \dot w &= R(v,w) \\
    v &= u_{\rm int}-u_{\rm ext}
\end{align*}
in $\Omega$, subject to boundary conditions, here 
\[
n^T\sigma_{\rm int}n^T\nabla u_{\rm int}=\sigma_{\rm ext}\nabla u_{\rm ext} = 0 \quad\text{on $\partial\Omega$}
\]  
and the ionic model dynamics $\dot w = R(u_{\rm int}-u_{\rm ext},w)$ in $\Omega$. Here, $\beta$ is the membrane area per volume, and the anisotropic diffusion tensors $\sigma_{\rm int}$ and $\sigma_{\rm ext}$ correspond to the conductivities $\sigma_i$, $i>0$, and $\sigma_0$, respectively, but are modified due to the geometric structure of the actual myocytes and their arrangement.

Besides the EMI model~\eqref{eq: cell-by-cell model}, we will rather focus on the monodomain equations, which can be derived from the bidomain system under the wrong but useful assumption of linearly dependent diffusion tensors $\sigma_{\rm int} = \lambda\sigma_{\rm ext}$, $\lambda\in\R$. The monodomain model is formulated directly in terms of the transmembrane voltage $v = u_{\rm int}-u_{\rm ext} \in H^1(\Omega)$ and gating variables $w\in L^2(\Omega)$ as
\begin{align}
    \beta\left(C_m \dot v + I^{\rm ion}(v,w)\right) &= \divergence(\sigma_m \nabla v) \label{eq:monodomain} \\
    \dot w &= R(v,w),
\end{align}
subject to boundary conditions $n^T\sigma_m\nabla v = 0$.

\paragraph{Weak formulation and finite element discretization}

Multiplying~\eqref{eq:monodomain} by a test function $\phi\in H^1(\Omega)$ yields the standard weak formulation
\[
    \int_\Omega \beta(C_m \dot v + I^{\rm ion})\phi + \nabla v^T \sigma_m \nabla \phi \, dx = 0
\]
and its standard discretization with Lagrange finite elements the ordinary differential equation system
\begin{equation} \label{eq: abstract monodomain}
\begin{aligned}
    M_m \dot v_h &= - A_m v_h - b_m(v_h,w_h) \\
    \dot w_h &= R(v_h,w_h).
\end{aligned}
\end{equation}
In contrast to the discretized EMI system~\eqref{eq: abstract reaction diffusion} and~\eqref{eq: abstract gating}, here the two ODEs have the same dimension $N_h$. 

\section{Time discretization}
\label{sec:SDC}

The discretized EMI and monodomain models~\eqref{eq: abstract reaction diffusion}--\eqref{eq: abstract gating} and~\eqref{eq: abstract monodomain}, respectively, both assume the abstract form
\begin{equation}\label{eq:abstract ODE}
\begin{aligned}
  B \dot z &= F(z)
\end{aligned}
\end{equation}
of an index one differential algebraic equation for $z=[u,w]\in \R^{N_h+n_h}$ or $z=[v,w]\in\R^{2N_h}$, respectively. For the purpose of the present Sec.~\ref{sec:SDC}, a distinction between EMI and monodomain model is not necessary.

In a method of lines approach, \eqref{eq:abstract ODE} can be integrated by, e.g., any L-stable single step method. For the purpose of cheap algebraic adaptivity worked out in Sec.~\ref{sec:adaptivity} below, we consider spectral deferred correction (SDC) methods~\cite{dutt2000spectral} in this paper.

\subsection{Spectral deferred correction methods}

SDC methods can be interpreted as implicit Runge-Kutta schemes on their own, but are essentially stationary iterations for solving collocation systems. On a single time step, w.l.o.g. $[0,T]$, we define a collocation time grid $0 < \tau_1 < \dots < \tau_m = T$ with $m$ collocation points, and ask for a polynomial approximation $z^*\in\Pol_m$ of $z$ of order $m$ respecting the given initial values $z(0)=z_0$ and satisfying~\eqref{eq:abstract ODE} at the collocation points $\tau_i$. Inclusion of the end point $T$ into the collocation grid while omitting the start point $t=0$ guarantees L-stability of the collocation solution~\cite{HairerNorsettWanner}. The most prominent example of such collocation grids are the Radau-IIa points used here, for which the resulting fully implicit Runge-Kutta schemes achieve convergence order $2m-1$.

Letting $\tau_0 = 0$ and $z_i = z(\tau_i)$, equation~\eqref{eq:abstract ODE} can be written as equivalent Picard equation and approximated by the high-order, i.e. spectral, quadrature rules $S_i$ corresponding to the collocation grid:
\begin{align} 
B (z_{i+1} - z_i) 
&= \int_{\tau=\tau_i}^{\tau_{i+1}} F(z)\,d\tau, \quad i=0,\dots,m-1 \notag \\
&\approx \sum_{j=0}^m S_{ij} F(z_j). \label{eq: SDC Integartion}
\end{align}
Applying Newton's method to~\eqref{eq: SDC Integartion} yields for $0\le i<m$
\begin{align*}
B(\delta z_{i+1}^k & - \delta z_i^k) - \sum_{j=0}^n S_{ij} F'(z_j^k) \delta z_j^k  \\
&= -B (z_{i+1}^k - z_i^k) + \sum_{j=0}^n S_{ij} F(z_j^k) =: \Phi(z^k)_i
\quad 
\end{align*}
for the current iterate $z^k$ and its Newton correction $\delta z^k$, yielding the  new iterate $z^{k+1} = z^k+\delta z^k$. In this equation system, all collocation time points are coupled due to $S$ being dense, which necessitates the expensive solution of a linear system of size $Nm$. Replacing $S$ on the left hand side by a suitable lower triangular matrix $\hat S$ decouples the collocation time points, such that $\delta z^k_{i+1}$ can be computed sequentially for $i=0,\dots,m-1$ in the SDC iteration:
\begin{align} \label{eq:SDC}
B(\delta z_{i+1}^k & - \delta z_i^k) - \sum_{j=0}^{i+1} \hat S_{ij} F'(z_j^k) \delta z_j^k  = \Phi(z^k)_i.
\end{align}
Due to the sequential progression through the collocation time points, the SDC iterations are often called sweeps. Different approximate quadrature rules $\hat S$ are in use, in particular the classical right-looking rectangular rule corresponding to an implicit Euler method and resulting in a bidiagonal $\hat S$~\cite{dutt2000spectral}, or the dense lower triangular $\hat S$ resulting from the so-called LU trick~\cite{Weiser2015} corresponding to diagonally implicit Runge-Kutta methods. 

Even though there is limited general convergence theory for arbitrary SDC methods, for specific choices of $m$, $\tau_i$, and $\hat S$, and for specific problem classes, SDC methods are well-known to converge with reasonable contraction factors $\rho < 0.6$.

\subsection{Operator splitting as basic scheme}

One of the attractive features of SDC is their flexibility in the choice of $\hat S$, allowing the use of problem-adapted basic solvers for~\eqref{eq:SDC}. For solving monodomain and EMI models, first order implicit-explicit operator splitting methods are popular. Using these methods as a basic integration method defining $\hat S$ allows exploiting the problem structure for cheaper SDC sweeps at the expense of a negligible or minor increase of the contraction rate $\rho$. Integrating operator splitting for monodomain and EMI models with LU-trick based SDC leads to the time discretization used in the present paper. In each iteration $k=1,\dots$, the values $\delta u^k_i$ of the potential correction (or transmembrane correction in case of monodomain) and $\delta w^k_i$ of the gating variable correction are computed for each collocation point $i=1,\dots,m$:
\begin{equation}\label{eq:final-sdc}
\begin{aligned}
&M(\delta u_{i+1}^k  - \delta u_i^k) + \sum_{j=0}^{i+1} \hat S_{ij} (A+\partial_u b(u^k_j,w^k_j)) \delta u_j^k  \\
&\; = -M (u_{i+1}^k - u_i^k) - \sum_{j=0}^n S_{ij} (Au^k_j+b(u^k_j,w^k_j)) \\
& u^{k+1}_{i+1} = u^k_{i+1} + \delta u^k_{i+1} \\
&(\delta w_{i+1}^k - \delta w_i^k) + \sum_{j=0}^{i+1} \hat S_{ij}\partial_w R(u_j^{k+1},w_j^k) \\
&\;= - w_{i+1}^k - w^k_i - \sum_{j=0}^n S_{ij} R(u^{k+1}_j,w^k_j) \\
& w^{k+1}_{i+1} = w^k_{i+1} + \delta w^k_{i+1} .
\end{aligned}
\end{equation}

In order to facilitate an efficient construction of the system matrices $\tilde A = M + \hat S_{jj}(A+\partial_u b(u^k_j,w^k_j))$, the variable reaction contribution $\partial_u b(u^k_j,w^k_j)$ is lumped and applied only on the diagonal.

Termination of the SDC iteration is controlled by the condition 
\[
\sum_{i=1}^m \|\delta u^k_i\|_{\tilde A}^2 \le \left(\frac{1-\rho}{\rho} \TOL_{\rm SDC}\right)^2,
\]
i.e. a sufficiently small energy norm error as estimated by the geometric series for the linearly convergent iteration.



\section{Nested subdomain selection}
\label{sec:adaptivity}

The locality of the excitation propagation dynamics, and in particular of the SDC corrections, see Fig.~\ref{fig:essential-support} and Figs.~\ref{fig:cells3}--\ref{fig:cells4} below, calls for spatio-temporal adaptivity in order to reduce the computational effort. Due to the travelling front like behavior of the solution, the SDC corrections are essentially confined to a neighborhood of the depolarization and repolarization fronts, and hence to a small region of the computational domain, see Fig.~\ref{fig:essential-support}.  The main ansatz investigated here is therefore, to restrict the computation of later SDC corrections to these relevant regions, and neglect the other parts of the domain $\Omega$. This type of spatial adaptivity can be done on the algebraic level, and due to lacking mesh modifications is cheap enough to reduce the computational effort considerably. 

\begin{figure}
    \centering
    \includegraphics[width=0.8\linewidth]{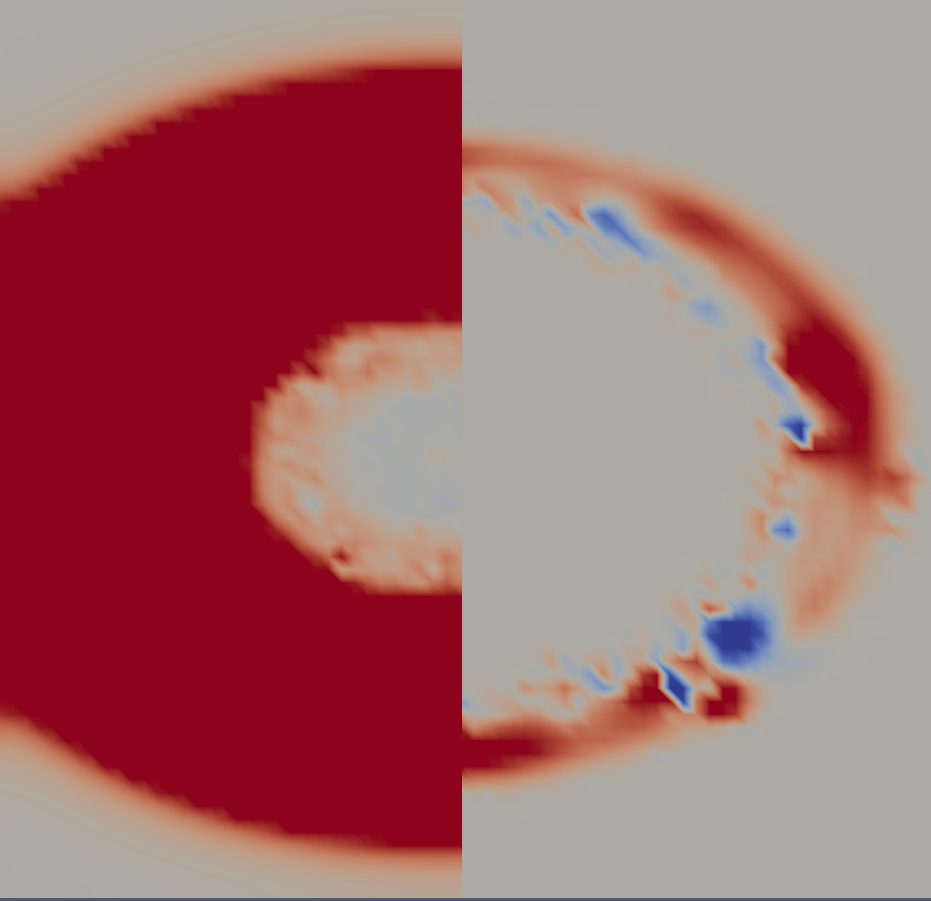}
    \caption{SDC corrections for a depolarization front in a 2D monodomain model spreading out from the center. Colour coded is the sign and magnitude of the correction $\delta v^k$. Everything in gray corresponds to a quantitatively negligible magnitude close to machine precision. \emph{Left:} first sweep. \emph{Right:} sweep 12. Not only is the essential support limited to a neighborhood of the depolarization front, it also shrinks with higher sweep number.}
    \label{fig:essential-support}
\end{figure}

The ansatz can also be interpreted as a kind of multirate integration: instead of using different time step sizes, different sweep counts and therefore convergence orders are used locally. In contrast to other predictor-corrector type multirate integration schemes~\cite{Rice1960,SavcencoHundsdorferVerwer2007}, the interleaving with the SDC iteration reduces the overhead of repeatedly integrating parts of the trajectory.

First we will derive an analytical a priori bound for the error introduced by this approach on the continuous level in Sec.~\ref{sec:aa-theo}. Guided by this result, we design an adaptive way of determining spatial subsets on the algebraic level and solve linear equation systems of reduced size in Sec.~\ref{sec:aa}.

\subsection{Continuous error bound} \label{sec:aa-theo}

In the monodomain case, a single SDC sweep~\eqref{eq:final-sdc} corresponding to fixed iteration index $k$ is a finite element and Runge-Kutta discretization of an inhomogeneous linear reaction-diffusion system for $(\delta v, \delta w)$.  Assuming the gating dynamics, i.e. the dynamics of $w$, to be slow compared to the transmembrane voltage dynamics, which is the case for several phenomenological ionic models, the equation for $\delta w$ can be neglected for the moment. We focus on the equation
\begin{equation} \label{eq:continuous-correction}
    \begin{aligned}
    \beta(C_m\dot{\delta v} &+ \partial_v I^{\rm ion}(v,w) \delta v) - \divergence(\sigma_m\nabla \delta v) \\
    &= \divergence(\sigma_m\nabla v) - \beta(C_m\dot v + I^{\rm ion}(v,w)) \\
    \delta v (0) &= 0
    \end{aligned}
\end{equation}
for the transmembrane voltage correction $\delta v$ and consider the restriction to a subdomain $\tilde\Omega\subset\Omega$ by applying homogeneous Dirichlet boundary conditions on $\partial\tilde\Omega$.
\newcommand{\T}{\mathopen]0,T\mathclose[}
\begin{theorem} \label{th:error-bound}
    Let $Q=\Omega \times \T$ and $\tilde Q=\tilde\Omega\times \T$. Assume $\delta v\in C^2(Q) \cap C(\bar Q)$ satisfies~\eqref{eq:continuous-correction} in $Q$ with homogeneous Neumann boundary conditions on $\partial\Omega\times \T$,  $\delta \tilde v\in C^2(\tilde Q) \cap C(\bar {\tilde Q})$ satisfies~\eqref{eq:continuous-correction} in $\tilde Q$ with homogeneous Dirichlet boundary conditions on $\partial\tilde\Omega\times\T$, and $I^{\rm ion}$ is continuously differentiable. Then the error $\epsilon = \delta v|_{\tilde\Omega} - \delta\tilde v$ is bounded by 
    \begin{equation*}
        \|\epsilon\|_{L^\infty(\tilde Q)} \le e^{\eta T} \|\delta v\|_{L^\infty(\partial\tilde\Omega \times \T)}
    \end{equation*}
    with $\eta = \max\big(0,-\min_{\tilde Q} \partial_v I^{\rm ion}\big)/(\beta C_m)$.
\end{theorem}
\begin{proof}
    Let $a = \partial_v I^{\rm ion}(v,w) / (\beta C_m)$ and $\kappa = \sigma_m/(\beta Cm)$. Then the error satisfies 
    \begin{align*}
        \dot\epsilon + a \epsilon - \divergence(\kappa\nabla\epsilon) &= 0 
        && \text{in $\tilde\Omega\times \T$} \\
        \epsilon &= 0 && t = 0 \\
        \epsilon &= \delta v && \text{on $\partial \tilde \Omega \times \T$}
    \end{align*}
     due to linearity of~\eqref{eq:continuous-correction}. Let $q = e^{\eta t}\epsilon$, which satisfies $\dot q + (a+\eta)q -\divergence(\kappa\nabla q)=0$ with $a+\eta\ge 0$ and Dirichlet boundary values $e^{\eta t}\delta v$. Applying the weak maximum principle~\cite{Evans2010} for parabolic equations yields 
     \begin{multline*}
         \|\epsilon\|_{L^\infty(\tilde Q)} \le \|q\|_{L^\infty(\tilde Q)} \\
         \le \|q\|_{L^\infty(\partial\tilde\Omega \times \T)} \le e^{\eta T}\|\delta v\|_{L^\infty(\partial\tilde\Omega \times \T)}
     \end{multline*}
     and therefore the claim.
\end{proof}
Theorem~\ref{th:error-bound} guarantees that the error introduced by restricting the computation of an approximate SDC correction $\delta\tilde v$ to $\tilde \Omega \subset \Omega$ and using homogeneous Dirichlet boundary conditions is bounded by $\|\delta v\|_{L^\infty(\bar\Omega\backslash\tilde\Omega)}$. 

Next we derive a bound on the error accumulation within the SDC iteration.

\begin{theorem}\label{th:r-bound}
    Assume the SDC iteration converges linearly, i.e. there is some $c\in\R$ such that $\|\delta v^k\|_{L^\infty(\Omega)} \le c\rho^k$. If the iteration count $r$ satisfies
    \begin{equation}\label{eq:r-condition}
    r \ge \frac{\log((1-\rho)\TOL/c)}{\log\rho} 
    \end{equation}
    and the approximate corrections $\delta \tilde v^k$ satisfy
    \[
    \|\delta \tilde v^k - \delta v^k\|_{L^\infty(\Omega)} \le \frac{1-\rho}{r+1}\TOL, \quad k=0,\dots,r,
    \]
    then the final approximation error meets the accuracy reqirement
    \begin{equation*}
        \|\sum_{k=0}^r \delta\tilde v - v \|_{L^\infty(\Omega)} \le \TOL
    \end{equation*}
    with at most one sweep more than the exact SDC iteration may take.
\end{theorem}
\begin{proof}
    Let $a=(1-\rho)/(r+1)$. By the triangle inequality, the final error is bounded by
    \begin{align*}
        \|\sum_{k=0}^r \delta\tilde v^k - v \|_{L^\infty(\Omega)}
        & \le \|\sum_{k=0}^r \delta\tilde v^k - \sum_{k=0}^r \delta v^k\|_{L^\infty(\Omega)} \\
        &\quad + \|\sum_{k=0}^r \delta v^k-  v \|_{L^\infty(\Omega)} \\
        &\le (r+1)a\TOL + c\frac{\rho^{r+1}}{1-\rho}.
    \end{align*}
    Due to~\eqref{eq:r-condition}, the last term equals $\rho\TOL$, such that we obtain
    \[
         \|\sum_{k=0}^r \delta\tilde v^k - v \|_{L^\infty(\Omega)}
         \le ((1+r)a+\rho)\TOL = \TOL.
    \]
    For the exact SDC iteration, $c\rho^{r+1}/(1-\rho)$ must not exceed $\TOL$ instead of $\rho\TOL$, such that one iteration less would be sufficient.
\end{proof}
Jointly, Theorems~\ref{th:error-bound} and~\eqref{th:r-bound} suggest, that the computational domain for $\delta\tilde v^k$ should be chosen as 
\begin{equation} \label{eq:subset}
\tilde\Omega_k = \{ x\in\Omega \mid |\delta v^k(x)| \ge \TOL_{\rm drop} \}
\end{equation}
with the drop tolerance 
\begin{equation} \label{eq:drop-tolerance}
    \TOL_{\rm drop} = e^{-\eta T} \frac{1-\rho}{r+1}\TOL
\end{equation}
being linear in the SDC iteration tolerance $\TOL$.

\subsection{Algebraic adaptivity} \label{sec:aa}

The subdomain selection suggested by~\eqref{eq:subset} can be realized approximately and efficiently on the algebraic level. When using Lagrangian finite elements associated to nodes $x_i \in \Omega$, instead of a subdomain $\tilde\Omega_{k+1}$ we define a subset of dofs, or indices,
\[
    I_{k+1} = \{ 1 \le i \le N_h \mid |\delta v^k(x_i)| \ge \TOL_{\rm drop} \}.
\]
Of course, the exact SDC correction $\delta v^k$ is unavailable except for the very first sweep, and thus we evaluate the approximate correction $\delta \tilde v^k$ computed on the index set $I_k$. Consequently, the index sets are nested, i.e. $I_{k+1}\subset I_k$, and lead to decreasing computational effort in the course of the iteration.

The value of $\TOL_{\rm drop}$ can in principle be computed explicitly from~\eqref{eq:drop-tolerance} based on estimates for SDC contraction $\rho$ and its error constant $c$.  Both can be estimated, even pointwisely, from monitoring the SDC convergence, which requires to perform at least two full sweeps and thus limits the possible speedup. Alternatively, since $\rho$ is relatively well known for reaction-diffusion systems~\cite{Weiser2015}, a reasonable value can be used and only $c$ estimated from the first sweep.

However, the resulting value for $\TOL_{\rm drop}$ may be too pessimistic, i.e. too small, for actual computation due to the worst case estimates in Theorems~\ref{th:error-bound} and~\ref{th:r-bound}. In the numerical experiments in Sec.~\ref{sec:experiments} below, we thus investigate a drop tolerance chosen according to $\TOL_{\rm drop} = \alpha\TOL$ with an empirically determined value of $\alpha$.

In each sweep, the equation systems~\eqref{eq:final-sdc} for $\delta v^k_{i+1}$ or $\delta u^k_{i+1}$ need to be solved. On the reduced set $I_k$ of dofs considered for sweep $k$, the submatrices $M_k$ and $A_k$ of $M$ and $A$, respectively, are needed. They can be obtained by an inexpensive submatrix extraction procedure. As the index sets $I_k$ are nested, this is a cheap progressive process.

One drawback of the submatrix selection is that preconditioners or factorizations need to be recomputed for every sweep, making expensive preconditioners and direct solvers less efficient. Since due to fast ionic currents the cardiac electrophysiology is usually only mildly stiff for reasonable spatial and temporal resolutions, we employ a conjugate gradient method with Jacobi preconditioner, which does not require any setup.

\section{Numerical experiments}\label{sec:experiments}

Here, we investigate the effectivity and efficiency of the algebraic adaptivity worked out in Sec.~\ref{sec:adaptivity} above on monodomain and EMI models both in 2D and 3D geometries. First we describe the experimental setup and then report on the adaptivity impact.

\subsection{Experimental setup}
\subsubsection{Common properties of all experiments}
For all numerical experiments we use the phenomenological Aliev-Panfilov model~\cite{ALIEV1996293} with $w\in\R$ for simplicity. The ion current is given by
\[
 I^{\rm ion}(v,w) = g_a v (v-a)(v-1) + vw,
\]
with the gating dynamics
\[
R(v,w) = \frac{1}{4} \left(\epsilon_1+\frac{\mu_1 w}{v+\mu_2}\right)(-w-g_s v(v-a-1)).
\]
The parameter values are $a = 0.1$, $\epsilon_1 = 0.01$, $g_a = 8.0$, $g_s = 8.0$ , $\mu_1 = 0.07$, $\mu_2 = 0.3$, such that the (dimensionless) transmembrane voltage covers the range $[0,1]$.

For the EMI models, the transmembrane current across the gap junctions in the intercalated discs connecting two myocytes is linear in the transmembrane voltage, and given by
\[
    I^{\rm ion}(v) = \frac{v}{R_g}.
\]
Spatial discretization is by linear finite elements, time discretization by SDC with LU trick on a RadauIIa grid with three collocation points, allowing up to order five convergence. The SDC tolerance $\TOL$ has been selected for all examples such that the truncation error is of the same magnitude as the time discretization error of the Radau collocation scheme. The initial guess for the SDC contraction factor is $\rho=0.05$.

The parameters of the monodomain and EMI models used are specified in Tab.~\ref{tab:parameters}. 

\begin{table}[h!]
	\label{tab:cell by cell parameters}
	\centering
	\begin{tabular}{@{\extracolsep{1pt}}lll}
		\toprule[1.5pt]
		Parameter & Value & Model\\
		\midrule 
		$\sigma_{>0}$ &$\SI{0.3}{\siemens\per\meter}$ & EMI\\
  		$\sigma_0$ & \SI{2.0}{\siemens\per\meter} & EMI \\
    	$C_m$  & \SI{1e-4}{\farad\per\square\meter} & EMI/MD\\
        $ R_g $  &$ \SI{4.5e-4}{\ohm.\square\meter}$ & EMI\\
       	$\sigma_m$ & \SI{0.3}{\siemens\per\meter} & MD \\
	    $ \chi $  & \SI{1400}{\per\meter} & MD\\
		\bottomrule[1.5pt]
	\end{tabular}
	\caption{Parameter used in the EMI and monodomain models, based on~\cite{JaegerEdwardsMcCullochTveito2019}.
    }
    \label{tab:parameters}
\end{table}

The linear systems~\eqref{eq:final-sdc} are solved by a Jacobi-pre\-con\-ditioned conjugate gradient solver, terminated at an estimated energy error reduction of $10^{-3}$. All computations have been performed using the Kaskade~7 finite element toolbox~\cite{GoetschelSchielaWeiser2020} compiled with GCC~10.2 on Linux 4.19. 2D examples have been run on a PC equipped with Intel Core i7-9700T CPU and \qty{32}{GB} RAM, 3D examples on a Dell PowerEdge 370 compute server with Intel Xeon E5 CPU and \qty{512}{GB} RAM.

\subsubsection{2D monodomain}
Here we set $\Omega=\mathopen]0,1\mathclose[^2$ discretized with a uniform Cartesian grid of mesh width $\delta x = \qty{6}{\micro\meter}$ resulting in \num{33025} mesh vertices, and use a time step of $T = \qty{1}{\milli\second}$ for a total integration time of \qty{0.5}{\second}. Excitation is initiated by setting $v = 0.5$ on the early excited domain $\{ x\in\Omega \mid \|x\| \leq 0.1 \}$.

\subsubsection{2D EMI}

\begin{figure}
  \centering
  \includegraphics[width=1\linewidth]{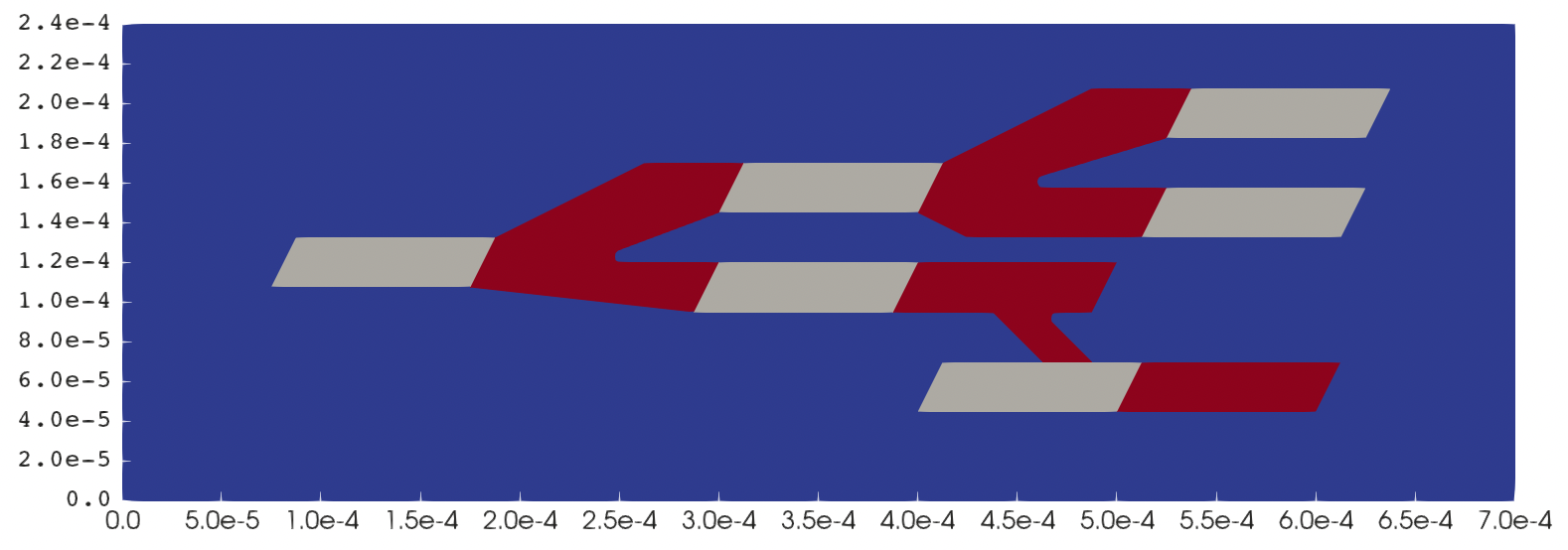}
  \caption{2D EMI setup. Myocytes coloured red/white, extracellular domain blue.}
  \label{fig:cells_2D}
\end{figure}

\begin{figure}
  \centering
  \includegraphics[width=0.9\linewidth]{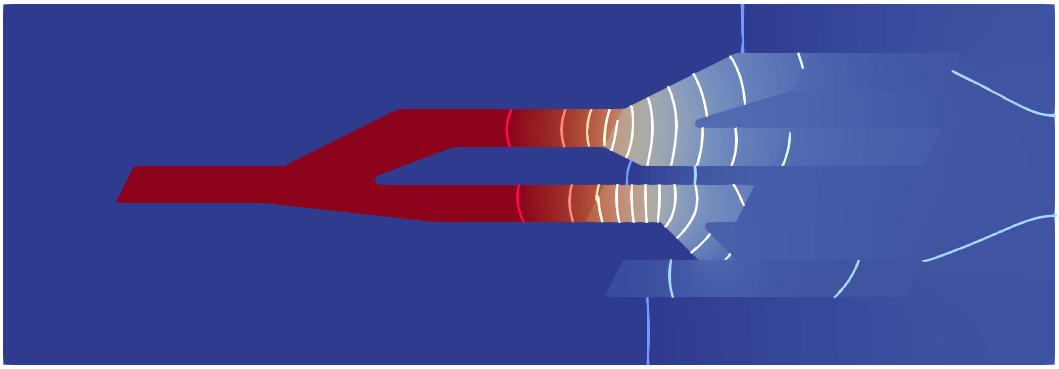}
  \caption{Potentials $u$ in the 2D EMI example after \SI{0.34}{\milli\second}. }
  \label{fig:cells2}
\end{figure}

As a simple cell-by-cell EMI setup, we consider the branched arrangement of myocytes in a domain of \qtyproduct{700 x 240}{\micro\meter} shown in Fig.~\ref{fig:cells_2D}. 
The domain is meshed with a triangle grid containing \num{129762} vertices.
Excitation is initiated by the initial value of the intracellular potential of the leftmost cell being set to the activated state. For the EMI models, we restrict the attention to short times and depolarization, and hence omit the gating dynamics by fixing $w=0$. The time step is $T=\qty{20}{\micro\second}$ for a total integration time of \qty{0.68}{\milli\second}. The resulting solution in terms of the intra- and extracellular potentials are shown in Fig.~\ref{fig:cells2}.

\subsubsection{3D monodomain}
\begin{figure}
  \centering
  \includegraphics[width=0.4\linewidth]{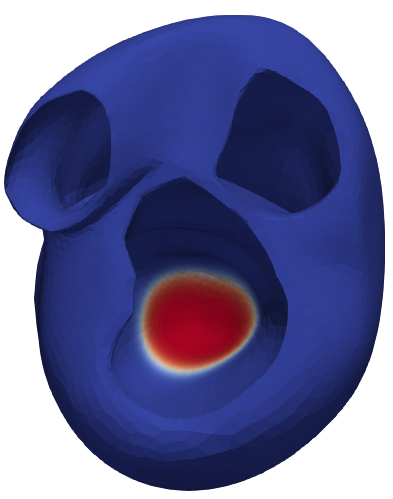}
  \includegraphics[width=0.5\linewidth]{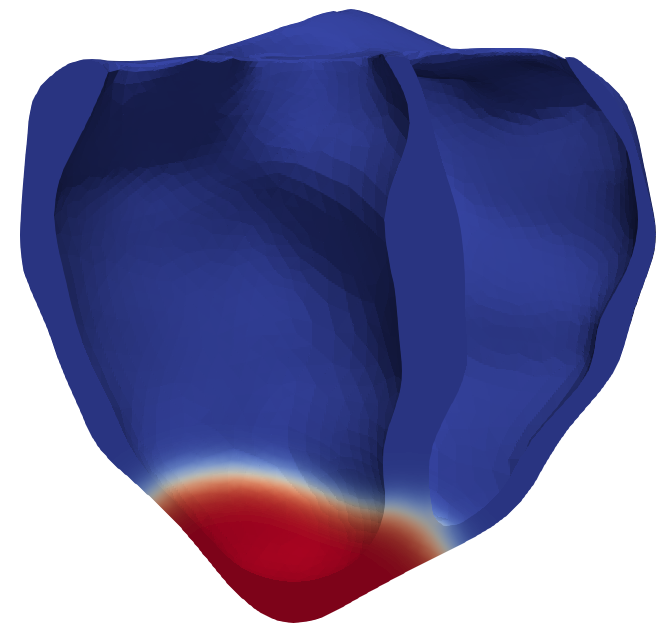}
  \caption{Transmembrane voltage in the 3D monodomain example after \qty{100}{\milli\second}.
  }
  \label{fig:cells3}
\end{figure}

For a more realistic 3D monodomian example, we use a generic ventricular geometry of rabbit heart dimension, i.e. of roughly \qtyproduct{1 x 1 x 1}{\centi\meter}, equipped with a tetrahedral mesh of \num{646166} vertices. Excitation is initiated at the apex. The resulting solution is shown in Fig.~\ref{fig:cells3}. 

\subsubsection{3D EMI}
\begin{figure}
 
  \centering
  \includegraphics[width=0.9\linewidth]{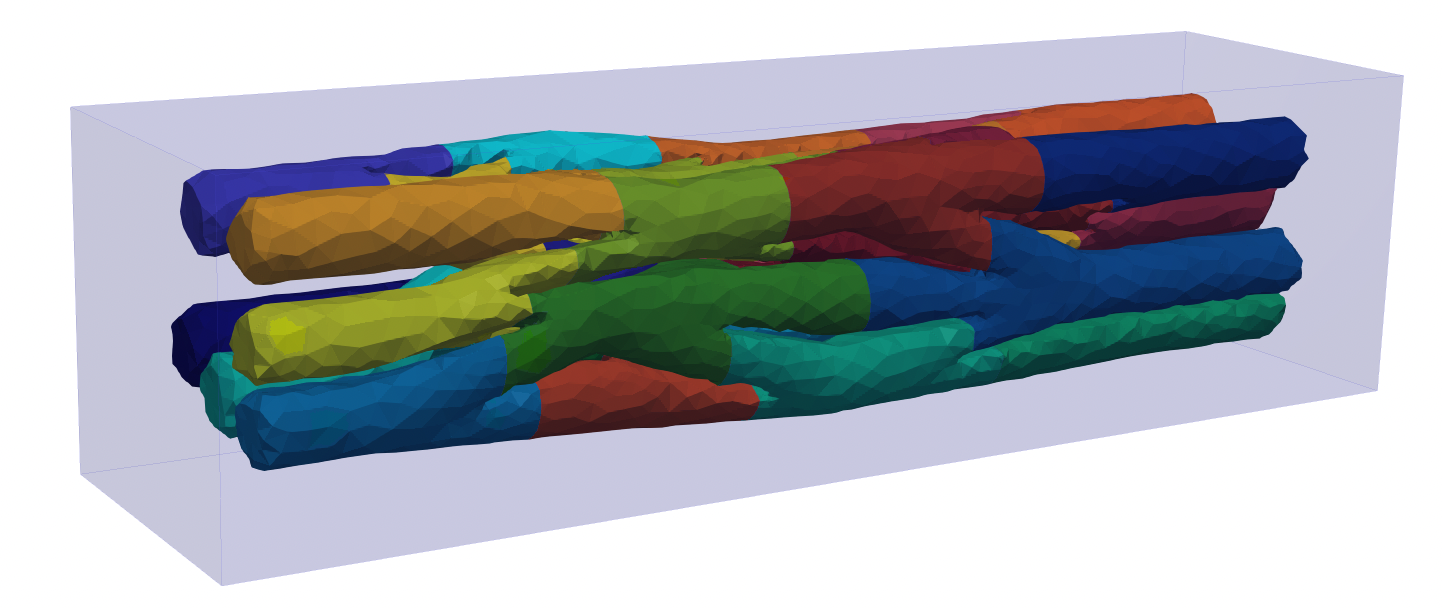}
  \caption{Geometry of the 3D EMI example with 43 myocytes surrounded by extracellular domain with dimension \qtyproduct{400 x 100 x 100}{\micro\meter} with \num{1923432} tetrahedra.}
  \label{fig:cell3D 43 cells}
\end{figure}

A complex 3D EMI setup consisting of 43 myocytes and extracellular domain is shown in Fig.~\ref{fig:cell3D 43 cells}, equipped with a tetrahedral mesh of \num{401120} vertices.
 \begin{figure}
  \centering
  \includegraphics[width=0.9\linewidth]{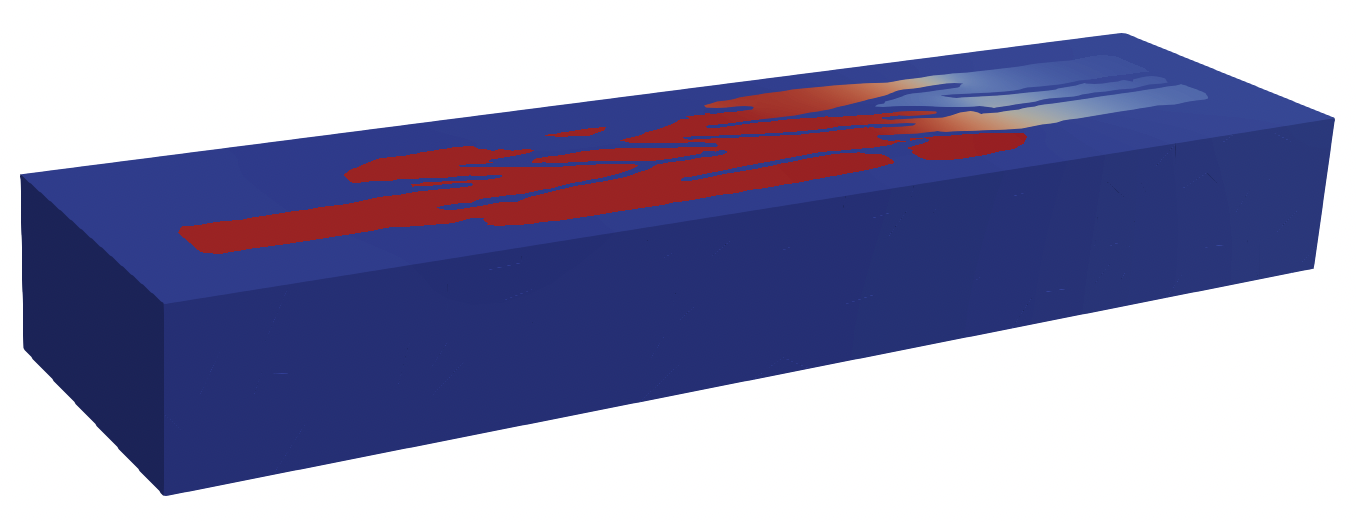}
  \caption{Solution of the 3D EMI example at $t = \qty{0.74}{\milli\second}$ (cross section)}
  \label{fig:cell3D 43 cells solution}
\end{figure}
Excitation is initiated by setting the front left bottom myocyte to the activated state.


\subsection{Algebraic adaptivity}

\begin{figure}
  \centering
  \includegraphics[width=1.0\linewidth]{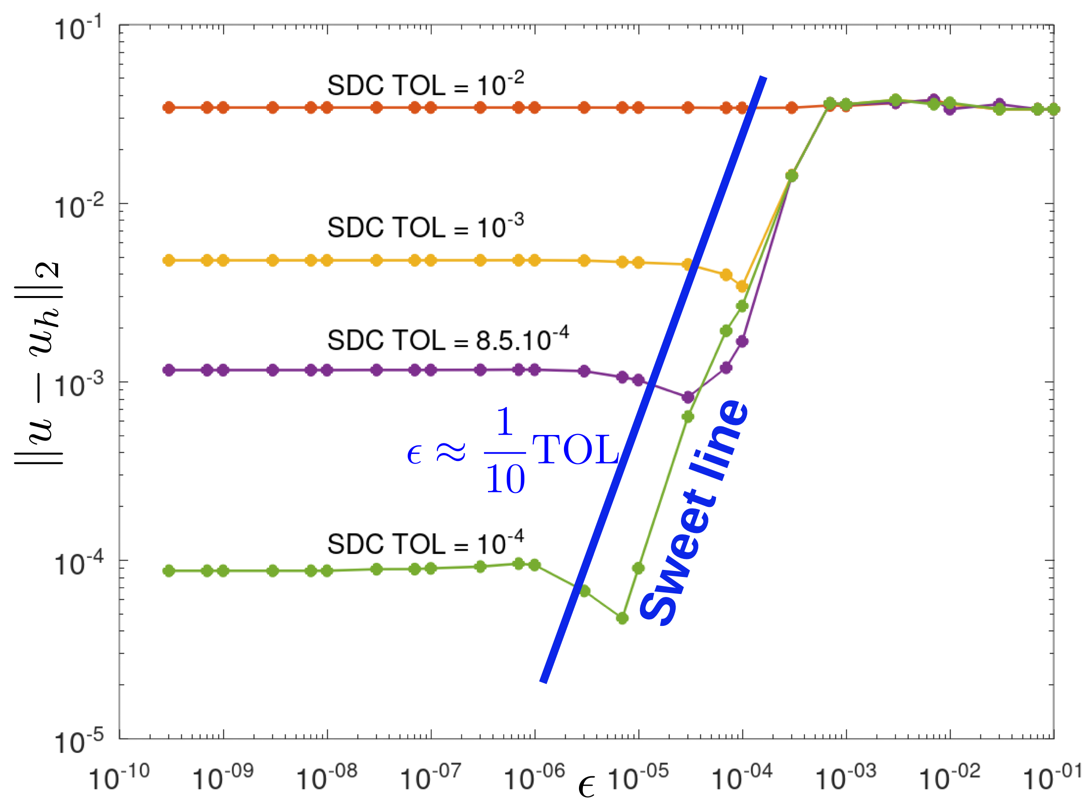}
  \caption{Simulation error in the 2D for EMI model at $t=\qty{2.5}{\milli\second}$ versus drop tolerance for different values of the SDC truncation tolerance $\TOL$. 
}
 \label{fig:cell_sweetline}
\end{figure}

First we investigate the impact of the drop tolerance $\TOL_{\rm drop}$ on the final accuracy in the 2D EMI example by comparing a wide range of combinations of SDC iteration tolerance $\TOL$ and drop tolerance, see Fig.~\ref{fig:cell_sweetline}.
As suggested by the theoretical result~\eqref{eq:drop-tolerance}, there is a threshold for $\TOL_{\rm drop}$ depending linearly on $\TOL$, such that the overall accuracy is barely affected by smaller drop tolerances, but dominated by the index subset selection error for higher drop tolerances. The numerical results indicate that the choice $\TOL_{\rm drop}=\frac{1}{10}\TOL$ is close to optimal, i.e. the maximum drop tolerance that does not affect solution accuracy.

\begin{figure}
  \centering
  \includegraphics[width=1.0\linewidth]{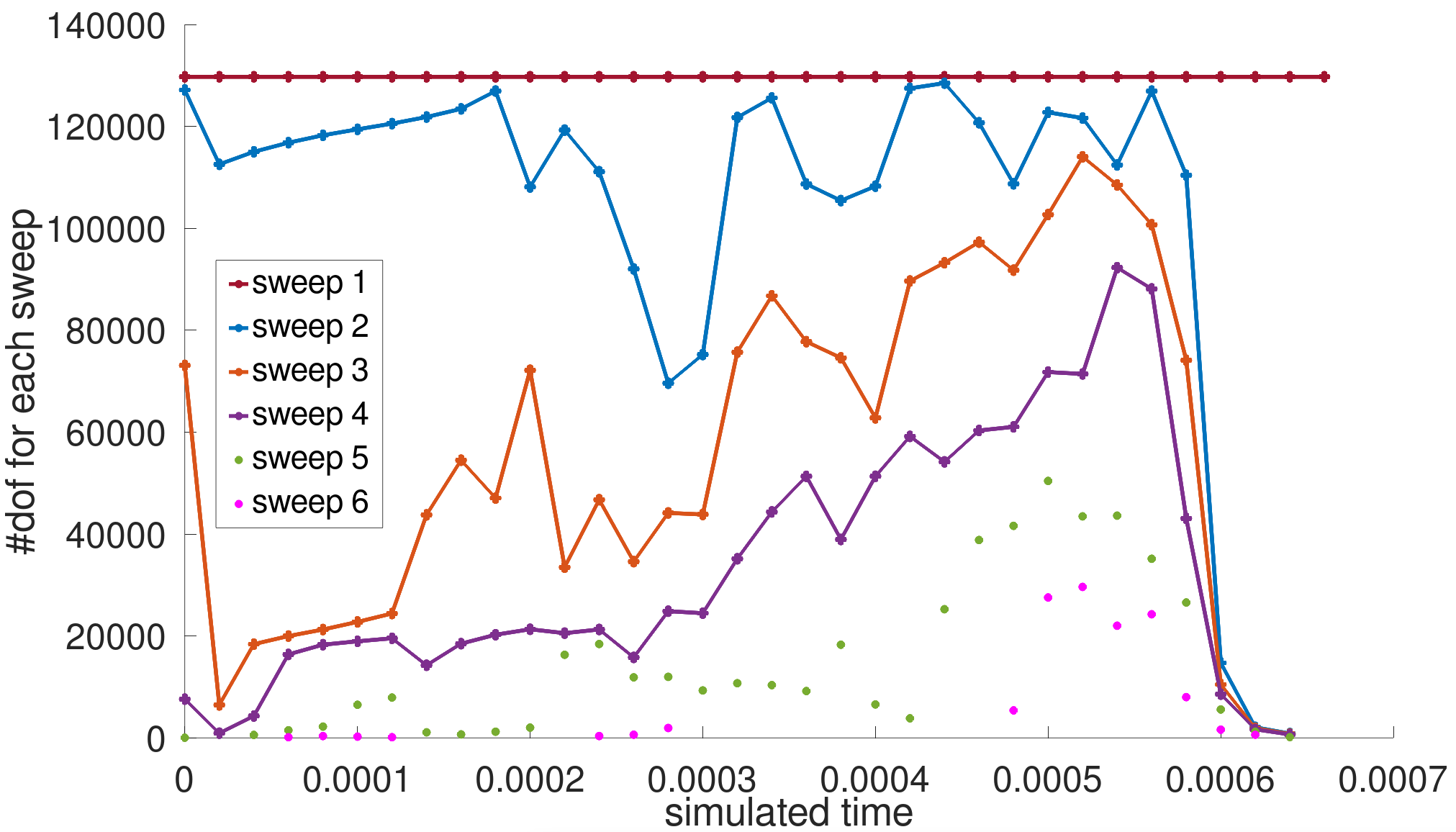}
  \caption{Number of active dofs in each sweep versus \qty{5}{\milli\second} simulated time $t$ in the 2D EMI example. SDC tolerance $\TOL=10^{-4}$, drop tolerance $\TOL_{\rm drop}=10^{-5}$ (on the sweet line). The  first sweep always uses all \num{129762} dofs. 
}
  \label{fig:cell_sweetline_time}
\end{figure}

The corresponding reduction of the linear equation size to solve during the sweeps is shown in Fig.~\ref{fig:cell_sweetline_time} versus the simulated time.

 \begin{figure}
  \centering
  \includegraphics[width=1.0\linewidth]{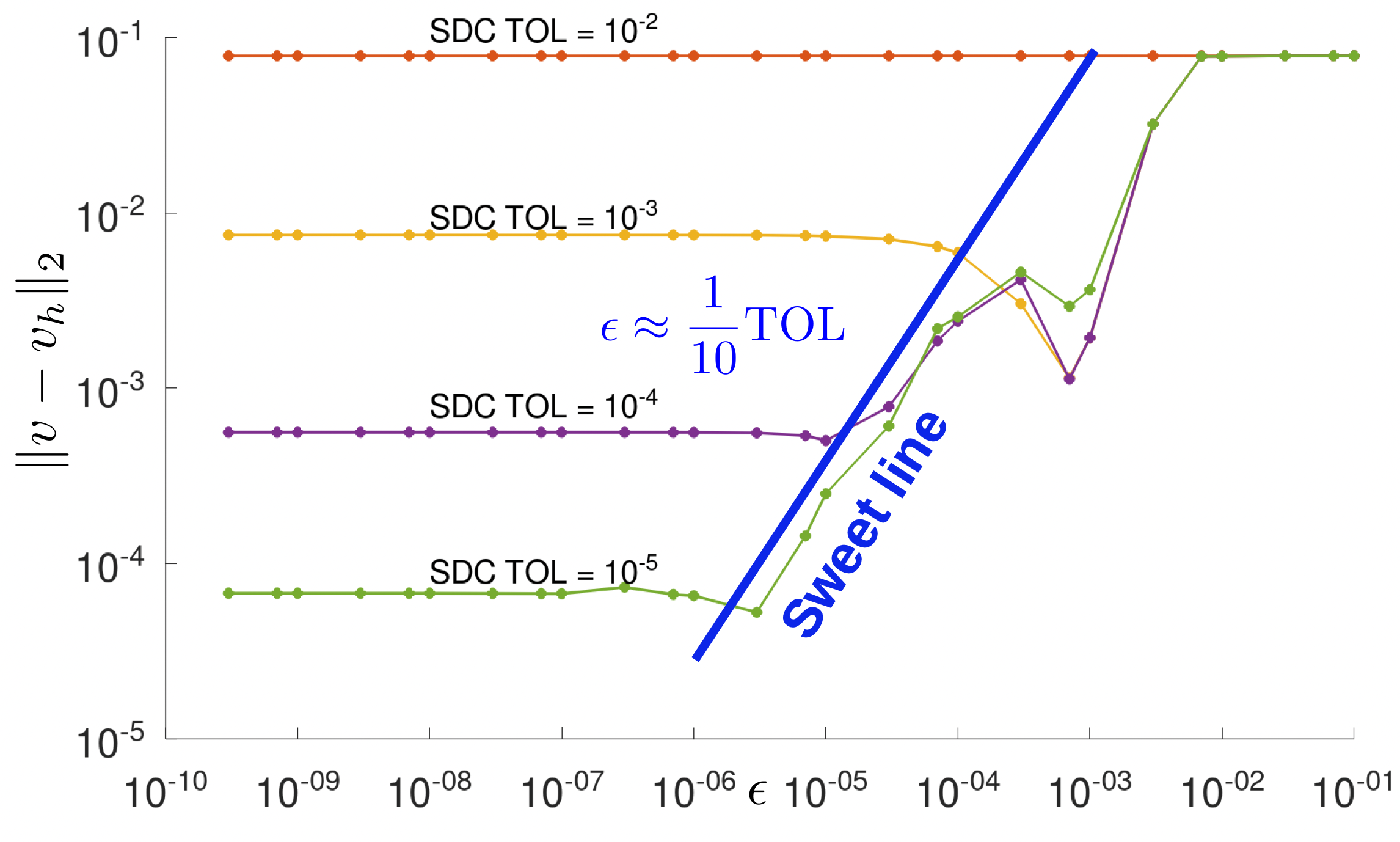}
  \caption{Simulation error in the 2D monodomain model at $t=\qty{100}{\milli\second}$ versus drop tolerance for different values of the SDC truncation tolerance $\TOL$. 
}
 \label{fig:mono_sweetline}
\end{figure}

The number of degrees of freedom retained in the later sweeps depends directly on the length of the depolarization front. This is clearly seen in the increase of active degrees of freedom in sweeps 3 to 6 over the integration time due to the arrangement of myocytes branching out, see Fig.~\ref{fig:cells_2D}, which leads to a longer front. At \qty{6}{\milli\second}, the depolarization front has traversed the whole domain, leading to a slow dynamics that can be captured perfectly well with just the first sweep.

In any case, the reduction of the number of active degrees of freedom, and hence of the problem size, is considerable: while the second sweep still contains most dofs, their number is halved in the third sweep and reduced even more for the subsequent sweeps.

\begin{figure}
  \centering
  \includegraphics[width=0.4\linewidth]{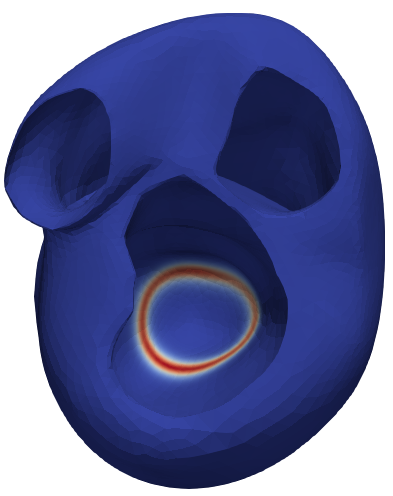}
  \includegraphics[width=0.5\linewidth]{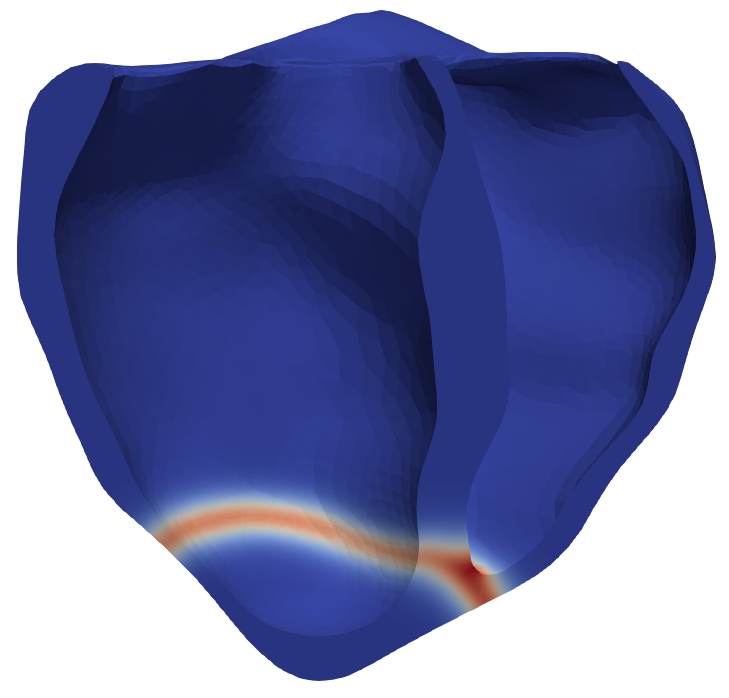}
  \includegraphics[width=0.4\linewidth]{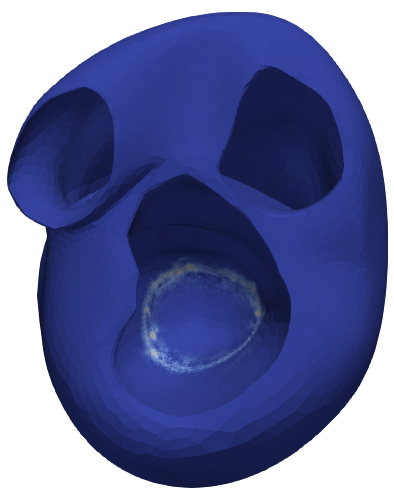}
  \includegraphics[width=0.5\linewidth]{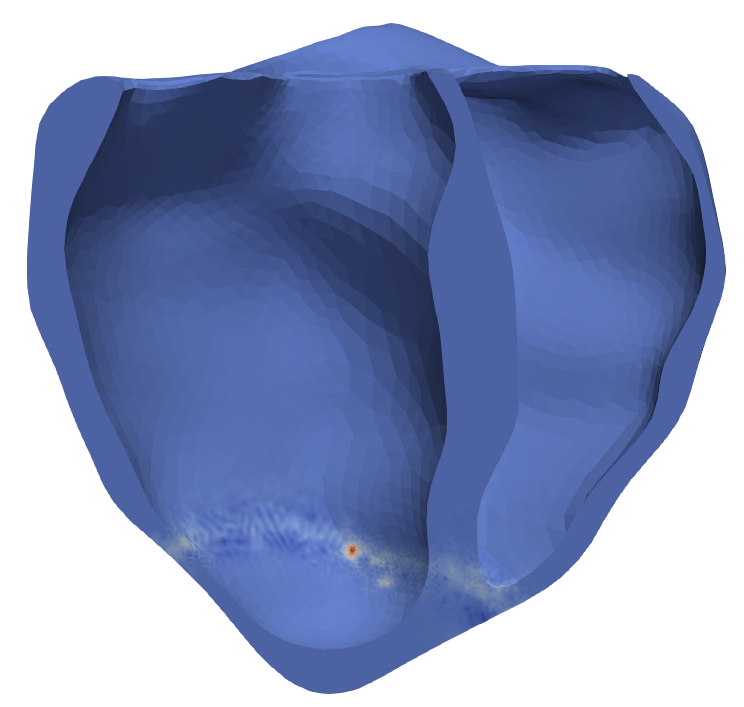}
  \includegraphics[width=0.4\linewidth]{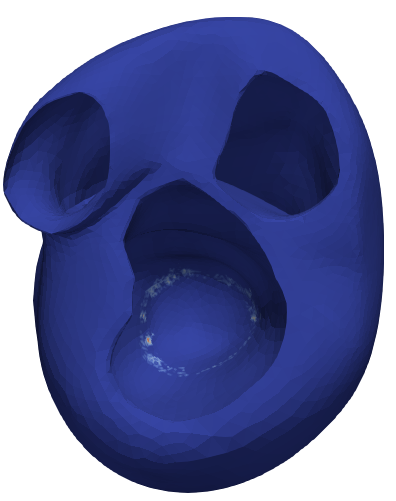}
  \includegraphics[width=0.5\linewidth]{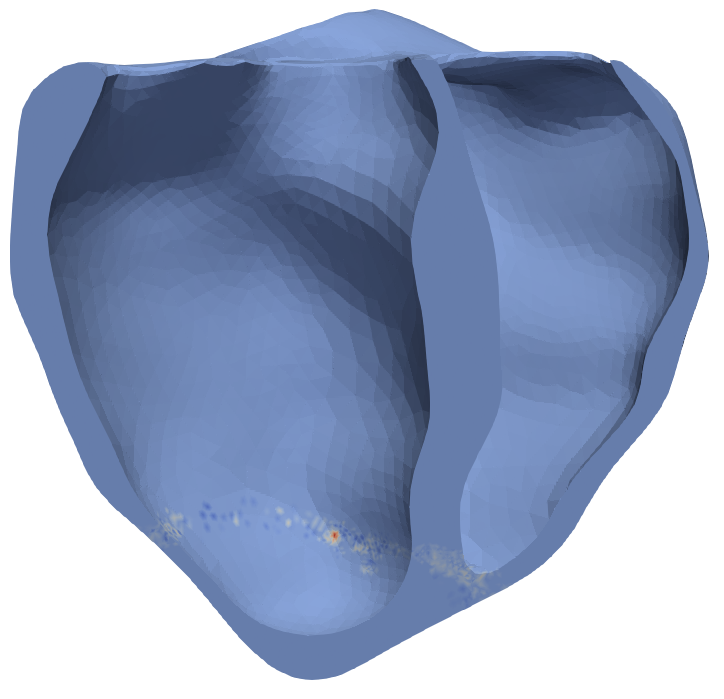}
  \caption{First to third SDC correction at final collocation point in the 3D monodomain example at $t=\qty{100}{\milli\second}$ (top to bottom). For the second sweep, \num{3283} dofs are selected out of \num{646166}, for the third sweep \num{425} dofs, and for the fourth sweep (not shown) \num{58}.}
  \label{fig:cells4}
\end{figure}

\begin{figure}
  \centering
  \includegraphics[width=1.0\linewidth]{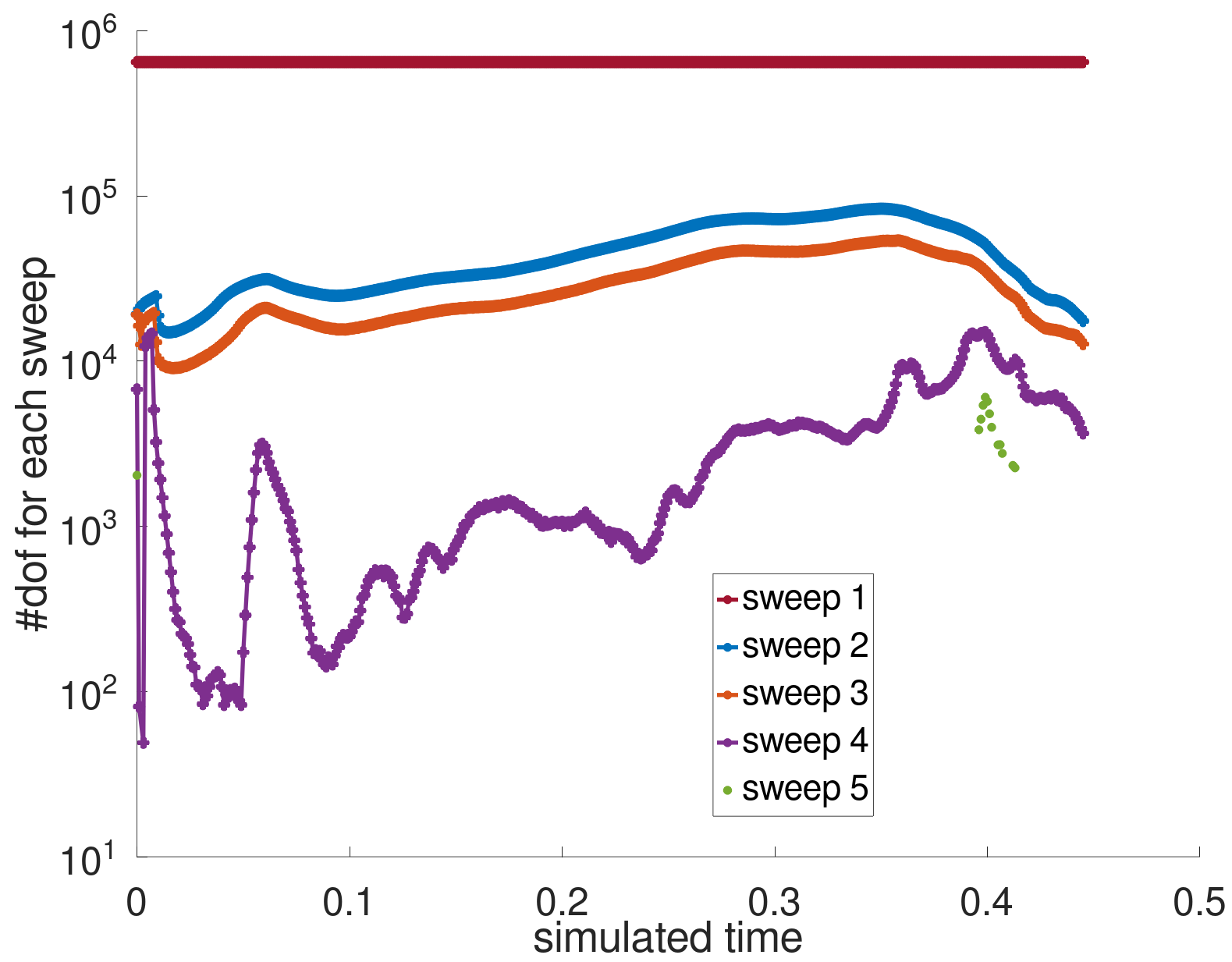}
  \caption{Number of active dofs in each sweep versus \qty{450}{\milli\second} simulated time $t$ in the 3D monodomain example. SDC tolerance $\TOL=10^{-4}$, drop tolerance $\TOL_{\rm drop}=10^{-5}$ (on the sweet line). The  first sweep always uses all \num{646166} dofs.}
  \label{fig:mono_sweetline_time}
\end{figure}

The same effect can be observed in the 3D monodomain example, see Fig.~\ref{fig:cells4}. 
The SDC corrections are essentially restricted to the fringe of the expanding activated region, which, in particular in the beginning of the excitation, is a relatively small part of the domain. Consequently, the index subset selection for later sweeps considerably reduces the linear equation systems' size, as shown in Fig.~\ref{fig:mono_sweetline_time}.

For all experiments, the simulation wall clock time has been measured with and without adaptivity. The results are given in Tab.~\ref{tab:speedups}, and show a significant speedup of the integration achieved by algebraic adaptivity without sacrificing accuracy. 

\begin{table}
\begin{center}
\begin{tabular}{ lcc } 
\toprule
  & monodomain & EMI \\
\midrule
   2D  & 2.13  &   3.336 \\ 
   3D  & 3.29   & 4.344 \\
   \bottomrule
\end{tabular}
\end{center}
   \caption{Computational speedup due to algebraic adaptivity with drop tolerance selected on the (empirical) sweet line. 
   }\label{tab:speedups}
\end{table}

\section{Discussion} \label{sec:discussion}

As illustrated by the numerical results in Sec.~\ref{sec:experiments} above, the proposed algebraic adaptivity is effective in reducing the number of degrees of freedom used in intermediate computations. In contrast to common mesh adaptivity approaches, the overhead is negligible, such that the smaller size of systems to solve translates directly in a significant reduction of simulation time compared to the baseline.

It is also clear, and observed in numerical experiments, that a larger spatial domain size, and thus a smaller fraction of the domain covered by the depolarization front, translates into higher speedup of the proposed adaptive approach.

The methodology has some restrictions, though. Due to the frequently changing linear equation systems to be solved, preconditioners with expensive setup cannot be used except for the very first SDC sweep. Due to locality of SDC corrections and moderate stiffness of monodomain equations, this appears not to be very restrictive. For EMI and bidomain models, block Jacobi preconditioners and a block selection procedure might be a reasonable path to follow. 

Second, as all spatial adaptivity approaches, the presented algebraic adaptivity incurs load imbalance in distributed simulations of large scale problems, leading to significantly reduced speedup. Dynamic load balancing can compensate that, but incurs some overhead, that itself reduces speedup. Nevertheless, algebraic adaptivity can be advantageous also in statically distributed simulations by reducing the energy consumption of large scale simulations.

Third, the speedup is currently limited by the number of SDC sweeps taken in the non-adaptive baseline method, since the first sweep is always performed on the whole domain. Since for reasonable collocation discretizations the number of SDC iterations can be expected to lie between 3 and 10, this limits the possible speedup considerably. Several directions for increasing the speedup further could be considered in future: increasing the collocation order and time step size, with a trade-off between SDC convergence rate, iteration count, and time step size, or ladder methods~\cite{layton65087implications,minion2003semi}, a cascadic multigrid variant using coarser collocation grids on the first sweeps, or even some heuristic a priori degree of freedom subset selection for the first sweep.

Finally, we would like to point out that the baseline used for comparison, i.e., the non-adaptive SDC operator splitting method, can but need not be the fastest algorithm for a given problem, such that the observed speedup has to be interpreted with due care~\cite{Goetschel_2021}.

\section*{Conclusions}
\label{sec:conclusions}

The combination of spectral deferred correction method and algebraic adaptivity by progressive vertex subset selection for later SDC sweeps is effective in reducing the computational cost of electrophysiology simulations both in simple monodomain and in more complex EMI models. Compared to mesh adaptivity, the procedure has a negligible overhead, such that the reduction in the number of considered degrees of freedom translates directly into efficiency gains. Speedup factors between 2 and 4 compared to the non-adaptive baseline have been observed.

\section*{Acknowledgements}
The authors would like to express their gratitude to Mark Potse for providing the 3D cell-by-cell grid.


\bibliographystyle{spmpsci}      
\bibliography{refs}   

\end{document}